\documentclass{commat}

\author{Jian Cui, Peter Danchev}
\affiliation{
\address{Jian Cui -- Department of Mathematics, Anhui Normal University, Wuhu, Anhui 241002, China}
\email{cui368@ahnu.edu.cn}

\address{Peter Danchev -- Institute of Mathematics and Informatics, Bulgarian Academy of Sciences, "Acad. G. Bonchev" str., bl. 8, 1113 Sofia, Bulgaria}
\email{danchev@math.bas.bg; pvdanchev@yahoo.com}
}

\title[Strongly $\pi$-$*$-Regular Rings] {On Strongly $\pi$-Regular Rings with Involution}



\abstract{Recall that a ring $R$ is called \emph{strongly $\pi$-regular} if, for every $a\in R$, there is a positive integer $n$ such that $a^n\in a^{n+1}R\cap Ra^{n+1}$. In this paper we give a further study of the notion of a \emph{strongly $\pi$-$*$-regular ring}, which is the $*$-version of strongly $\pi$-regular rings and which was originally introduced by Cui-Wang in J. Korean Math. Soc. (2015). We also establish various properties of these rings and give several new characterizations in terms of (strong) $\pi$-regularity and involution. Our results also considerably extend recent ones in the subject due to Cui-Yin in Algebra Colloq. (2018) proved for $\pi$-$*$-regular rings and due to Cui-Danchev in J. Algebra Appl. (2020) proved for $*$-periodic rings.}

\VOLUME{31}
\YEAR{2023}
\NUMBER{1}
\firstpage{73}
\DOI{https://doi.org/10.46298/cm.10273}

\begin{paper}

\section{Introduction and Background}

Everywhere in the text of the present paper, all rings into consideration are assumed to be associative and containing the identity element $1$ which differs from the zero element $0$, and all subrings are unital (i.e., containing the same identity as that of the former ring). Our terminology and notations are mainly standard being in agreement with \cite{L}. Concretely, $U(R)$ denotes the group of all units in $R$, $Id(R)$ the set of all idempotents in $R$, $Nil(R)$ the set of all nilpotents in $R$, and $J(R)$ the Jacobson radical of $R$. The more specific notions will be given in detail. To explain what this might mean, we recall the notion of a $*$-ring. An involution of a ring $R$ is an operation $* : R \to R$ such that $(x + y)^* = x^* + y^*$, $(xy)^* = y^*x^*$ and ${(x^*)}^* = x$ for all $x, y \in R$. Evidently, the identity mapping $\mathrm{id}_R$ is an involution of $R$ if and only if $R$ is commutative. A ring $R$ with involution $*$ is called a $*$-ring. An element $e$ of a $*$-ring $R$ is called a projection if $e^2 = e = e^*$.

We continue with the following concept (see \cite{CW}).

\begin{definition}\label{0} A $*$-ring $R$ is said to be \emph{$*$-abelian} if all projections of $R$ are central.
\end{definition}

Clearly, abelian $*$-rings that are $*$-rings for which each idempotent is central are just $*$-abelian, whereas the converse implication is untrue. Nevertheless, in accordance with \cite[Lemma~2.1]{LZ}, a $*$-ring is abelian provided if every its idempotent is a projection. Contrastingly, in an abelian $*$-ring an arbitrary idempotent need not be a projection. In fact, if $R = S\times S$, where $S$ is a commutative ring, with involution $*$ given by $(a, b)^* = (b, a)$ for some $a,b\in R$, then $R$ is a $*$-ring that is abelian, but the idempotent $(1, 0)$ is obviously not a projection. Note that Lemma~\ref{3} listed below explains when an abelian $*$-ring has the property that every idempotent is a projection.

Likewise, if $R$ is a $*$-ring, then the quotient $R/J(R)$ is also a $*$-ring, which fact is pretty easy and so we leave it without proof. However, it will be used freely in the sequel.

Some of the important achievements on certain properties of $*$-rings in a historical aspect are these: In \cite{V} was studied the $*$-version of clean rings and some other related classes of rings. After that, in \cite{LZ} was investigated the strong variant of $*$-cleanness. Later, in \cite{CY} and \cite{CXZ} were examined the $*$-$\pi$-regularity and $*$-nil-cleanness, respectively. Recently, in \cite{CD} were discovered the structure of $*$-periodic rings as well.

Our motivation in writing up this article is to explore the strong version of $*$-$\pi$-regular rings by raising some new ideas and approaches in establishing the major characterization theorem, which describes the structure of these strongly $*$-$\pi$-regular rings. The organizational plan of our work is as follows: We first of all begin with some technical claims (namely, Lemma~\ref{1}, Lemma~\ref{2}, Lemma~\ref{3} and Lemma~\ref{4}) which are our key instruments. Using them, we next are able to prove the basic statement, namely Theorem~\ref{5}. The subsequent construction in Example~\ref{6} unambiguously illustrates that the things are, in general, not too easy than it could be anticipated. Indeed, it is shown there the curious fact that there exists a non-commutative $*$-abelian strongly $\pi$-regular ring that is definitely \emph{not} strongly $\pi$-$*$-regular. This allows us to detect another strategy, given in the technical Lemma~\ref{7}, in attacking the proof of the complete criterion, namely Proposition~\ref{8}, for an arbitrary $*$-ring to be strongly $\pi$-$*$-regular, which is stated in a rather more suitable and comfortable for applications form. For completeness of the exposition and reader's convenience, we close with a few commentaries in Remark~\ref{9} and also pose two still open Problem~\ref{10} and Problem~\ref{11}, respectively.

\section{Preliminary and Main Results}

We start here with a series of four technicalities, which to keep a record straight also appeared in \cite{CXZ}, but we state them with slightly modified and simplified proofs only for the sake of completeness and for the readers' convenience.

\begin{lemma}\label{1} Let $R$ be a $*$-ring and $e^2 = e \in R$. If $1 + {(e - e^*)}^*(e - e^*) \in U(R)$, then
there exists $p^2 = p^* = p \in R$ such that $eR = pR$ and $e - p \in Nil(R)$.
\end{lemma}

\begin{proof} Set $x := 1 + {(e - e^*)}^*(e - e^*)$. Then, $xe = ee^*e = ex$ and $x^* = x$, so that $xe^* = e^*x$
and $e = x^{-1}ee^*e$. Therefore, $x^{-1}ee^* = x^{-1}ee^*x^{-1}ee^*$. Putting $p := x^{-1}ee^*$, we calculate that $p = p^2 = p^*$. Since $(1-e)p = (1-e)x^{-1}ee^* = 0$, one has that $p = ep$. However, $pe = (x^{-1}ee^*)e = x^{-1}(ee^*e) = x^{-1}xe =e$, so we get $eR = pR$ and $(e - p)^2 = e + p - ep - pe = 0$. Finally, it is straightforward to see that $e - p \in Nil(R)$, as wanted.
\end{proof}

\begin{lemma}\label{2} Let $R$ be a $*$-ring and $e^2 = e \in R$.
\begin{enumerate}
\item If $e - e^* \in J(R)$, then there exists $p^2 = p^* = p \in R$ such that $eR = pR$ and $e - p \in J(R) \cap Nil(R)$.

\item If $e - e^* \in Nil(R)$, then there exists $p^2 = p^* = p \in R$ such that $eR = pR$ and $e - p \in Nil(R)$.
\end{enumerate}
\end{lemma}

\begin{proof}
\begin{enumerate}
\item Since $e - e^* \in J(R)$, one follows by a straightforward check that 
$$x = 1 + {(e - e^*)}^*(e - e^*) \in U(R).$$
 Now, Lemma~\ref{1} enables us that there exists $p^2 = p^* = p \in R$ such that $eR = pR$ and $e - p \in Nil(R)$. So, $e = pe$ and hence $p = ep = p^* = pe^*$, which in turn assures that
$$e - p = pe - pe^* = p(e - e^*) \in J(R).$$
Thus, $e - p \in J(R) \cap Nil(R)$, as pursued.

\item Notice that ${(e - e^*)}^*(e - e^*) = (e - e^*){(e - e^*)}^*$. Therefore, 
$$x = 1 + {(e - e^*)}^*(e - e^*) \in U(R)$$ as $e - e^* \in Nil(R)$. Thus, the asked assertion follows directly from Lemma~\ref{1}.
\qedhere
\end{enumerate}
\end{proof}

\begin{lemma}\label{3} The following conditions are equivalent for a $*$-ring $R$:

\begin{enumerate}
\item For each $e^2 = e \in R$, $e - e^* \in J(R)$ and $R$ is $*$-abelian.\label{l3.1}

\item For each $e^2 = e \in R$, $e - e^* \in Nil(R)$ and $R$ is $*$-abelian.\label{l3.2}

\item Every idempotent of $R$ is a projection.\label{l3.3}

\end{enumerate}
\end{lemma}

\begin{proof} By virtue of \cite[Lemma 2.1]{LZ}, the ring $R$ is necessarily abelian if all idempotents of $R$ are projections. So, both the implications \ref{l3.3} $\Rightarrow$ \ref{l3.1}  and \ref{l3.3}  $\Rightarrow$ \ref{l3.2}  are clear. We next will verify only the implication \ref{l3.1}  $\Rightarrow$ \ref{l3.3} , because the remaining one \ref{l3.2}  $\Rightarrow$ \ref{l3.3}  can be shown similarly.
In fact, in view of Lemma~\ref{2}, there exists a projection $p \in R$ such that $eR = pR$. Since $R$ is $*$-abelian, we easily obtain that $e = pe = ep = p$, as required.
\end{proof}

\begin{lemma}\label{4} The following points are equivalent for a $*$-ring $R$:

\begin{enumerate}
\item For each $e^2 = e \in R$, $e - e^* \in J(R)$ and $ee^* = e^*e$.\label{l4.1}

\item For each $e^2 = e \in R$, $e - e^* \in Nil(R)$ and $ee^* = e^*e$.\label{l4.2}

\item Every idempotent of R is a projection.\label{l4.3}
\end{enumerate}
\end{lemma}

\begin{proof} \ref{l4.1}. $\Rightarrow$ \ref{l4.3}. Assume that \ref{l4.1} holds. Write $j := e - e^*$. Then $$e(1 - e^*) = j(1 - e^*) \in
J(R)$$ and $$e^*(1 - e) = -j(1 - e) \in J(R)\,.$$ As $ee^* = e^*e$, it readily follows that both $e(1 - e^*)$ and $e^*(1 - e)$ are also idempotents, which in turn implies that $e = ee^* = e^*e = e^*$. So every idempotent of $R$ is a projection, as stated.

The verification of \ref{l4.2}. $\Rightarrow$ \ref{l4.3}.~is analogous to that of \ref{l4.1}.  $\Rightarrow$ \ref{l4.3}., and the proofs of the implications \ref{l4.3}. $\Rightarrow$ \ref{l4.1}.~and \ref{l4.3}. $\Rightarrow$ \ref{l4.2}.~are both rather clear, so that we leave them to the interested reader for a direct check.
\end{proof}

Let us recall that a ring $R$ is called \emph{$\pi$-regular} (cf.~\cite{T}) if, for every $a \in R$, the relation $a^n \in a^nRa^n$ holds for some integer $n\geq 1$ and is called \emph{strongly $\pi$-regular} if, for every $a\in R$, the relation $a^n\in a^{n+1}R\cap Ra^{n+1}$ holds for some integer $n\geq 1$. It is well known that strongly $\pi$-regular rings are always $\pi$-regular, while the converse is not true. However, for abelian rings that are rings whose idempotents are central, the equivalence is valid (see, e.g.~\cite{B}).

We have now accumulated all the information needed to prove the following criterion for a ring to be strongly $\pi$-$*$-regular.

\begin{theorem}\label{5} Let $R$ be a $*$-ring. Then the following statements are equivalent:

\begin{enumerate}
\item $R$ is strongly $\pi$-$*$-regular.

\item $R$ is $\pi$-regular and exactly one of the following equivalent conditions holds:

\begin{enumerate}
\item for each $e^2 = e \in R$, $e - e^* \in J(R)$ and $R$ is $*$-abelian;

\item for each $e^2 = e \in R$, $e - e^* \in Nil(R)$ and $R$ is $*$-abelian;

\item for each $e^2 = e \in R$, $e - e^* \in J(R)$ and $ee^* = e^*e$;

\item for each $e^2 = e \in R$, $e - e^* \in Nil(R)$ and $ee^* = e^*e$.
\end{enumerate}

\item $R$ is $*$-abelian and, for each $a \in R$, $Ra^m = Ra^m{(a^m)}^*a^m$ for some natural number $m$.

\item For each $a \in R$, there exist $e^2 = e \in R$ and an integer $n \geq 1$ such that
\[ Ra^n = R{(a^n)}^* = Re. \]
\end{enumerate}
\end{theorem}

\begin{proof}
\begin{description}
\item[1. $\Rightarrow$ 2.] It is apparent that $R$ is $\pi$-regular. In view of \cite[Theorem 3.6 (2)]{CW} or \cite{CY}, every idempotent of $R$ is a projection. So, the result follows utilizing both Lemmas~\ref{3} and \ref{4}.

\item[2. $\Rightarrow$ 3.] Applying Lemmas~\ref{3} and \ref{4}, every idempotent of $R$ is a projection. Thus $R$ is abelian, and so it is $*$-abelian. As abelian $\pi$-regular rings are always strongly $\pi$-regular, by using \cite[Proposition 1]{N}, one has that $a^m = eu$, where $e \in Id(R)$ and $u \in U(R)$ such that $a, e, u$ commuting
with each other. That is the reason why $e = a^mu^{-1} = u^{-1}a^m$ and $e{(u^{-1})}^* = {(u^{-1})}^*e$, because $e^* = e$.  Thus, 
\[
a^m = ea^m = ee^*a^m = e{(u^{-1})}^*{(a^m)}^*a^m = {(u^{-1})}^*e{(a^m)}^*a^m = {(u^{-1})}^*u^{-1}a^m{(a^m)}^*a^m,
\]
which yields that $Ra^m =Ra^m{(a^m)}^*a^m$, as expected.

\item[3. $\Rightarrow$ 4.] For any $e\in Id(R)$, one has that $Re^* = Re^*ee^*\subseteq Ree^* \subseteq Re^*$. So, $Re^*=Ree^*$. Letting $e^* = t^*ee^*$ for some $t \in R$, we then deduce that $$e^*t = t^*ee^*t = {(e^*t)}^*e^*t.$$ Writing $p := e^*t$, we infer that $p = p^*p = p^*$ and $p^2 = p^*p = p$. So, it easily follows that $Rp = Rp^* = Rt^*e \subseteq Re$, whence $Re = R{(e^*)}^* = Ree^*t \subseteq Rp$. Thus, $Re = Rp$. It next follows that $e = ep$ and $pe = p$. Since $R$ is $*$-abelian, we conclude that $e = p$. Consequently, every idempotent of $R$ is a projection. Let $a \in R$. By hypothesis, we have that $$Ra^m = Ra^m{(a^m)}^*a^m \subseteq R{(a^m)}^*a^m \subseteq Ra^m\,.$$ Therefore, $Ra^m = R{(a^m)}^*a^m$. Letting now $a^m = x{(a^m)}^*a^m$ for some element $x \in R$, we then obtain that $$a^mx^*a^m = (x{(a^m)}^*a^m)x^*a^m = x{(x{(a^m)}^*a^m)}^*a^m = x{(a^m)}^*a^m = a^m.$$
Write $f := x^*a^m$. Then, one inspects that $f^2 = f = f^*$ and $Ra^m = Rf$, because $Ra^m = Ra^mf \subseteq Rf \subseteq Ra^m$. Since $R$ is $*$-abelian, we have $$R{(a^m)}^* = R{(a^mf)}^* = R{(a^m)}^*f \subseteq Rf = Ra^m\,.$$
 Further, if we take $g := a^mx^*$, then $g^2 = g = g^*$ and $$a^m = ga^m = a^mg = a^mg^* = a^mx{(a^m)}^*\,.$$
 So, we finally arrive at $Ra^m \subseteq R{(a^m)}^*$ and, therefore, $Ra^m = R{(a^m)}^*$, as desired.

\item[4. $\Rightarrow$ 1.] For any $f \in Id(R)$, one observes that $Rf = Rf^*$. Then, $f = ff^*$, and hence it elementarily follows that $f = {(ff^*)}^* = f^*$. Thus, the idempotents of $R$ are projections. According to \cite[Lemma 2.1]{LZ}, the ring $R$ is abelian. By assumption, one can write that $a^n = a^ne$ and $e = ra^n$ for some element $r \in R$. Furthermore, one checks that $a^n = a^nra^n$, and thus $R$ is an abelian $\pi$-regular ring. In view of \cite{CY}, one writes that $a = b + pv = b + vp$, where $b \in Nil(R)$ with $ab = ba$, $p^2 = p = p^*$ and $v \in U(R)$. This proves finally that $R$ is strongly $\pi$-$*$-regular, as promised.
\qedhere
\end{description}
\end{proof}

The next example is crucial.

\begin{example}\label{6} There exists a $*$-abelian (strongly) $\pi$-regular ring which is not strongly $\pi$-$*$-regular. For
instance, let $R = \mathbb{T}_2(\mathbb{Z}_2)$ be the $2 \times 2$ upper triangular matrix over the two element field $\mathbb{Z}_2$, and an involution $* : R \to R$ defined by
$$
\begin{pmatrix} a&b\\0&c\end{pmatrix} \hookrightarrow \begin{pmatrix} c&b\\0&a\end{pmatrix}.
$$
for some $a,b,c\in \mathbb{Z}_2$. As $R$ has only the trivial projections, it is pretty easy to see that $R$ is $*$-abelian. However, it is plainly observed that $R$ is (strongly) $\pi$-regular, but it is manifestly \emph{not} strongly $\pi$-$*$-regular since $R$ is surely \emph{not} abelian. This substantiates our claim after all.
\end{example}

The next technical claim is useful.

\begin{lemma}\label{7} Let $R$ be a $*$-ring. If $J(R)$ is nil, then every projection of the factor-ring $R/J(R)$ is lifted to a projection of $R$.
\end{lemma}

\begin{proof} Let $\overline{e} = {\overline{e}}^* = \overline{e}^2 \in R/J(R)$. Since all idempotents lift modulo $J(R)$ as $J(R)$ is supposed to be nil, we may without loss of generality assume that $e^2 = e \in R$. Since ${\overline{e}}^* = \overline{e^*}$, one concludes that $e-e^* \in J(R)$. Due to Lemma~\ref{2}, there exists $p^2 = p^* = p \in R$ such that $e - p \in J(R)$, as desired.
\end{proof}

The following necessary and sufficient condition can be somewhat treated as a supplement to the chief Theorem~\ref{5} proved above.

\begin{proposition}\label{8} Suppose $R$ is a $*$-ring. Then $R$ is strongly $\pi$-$*$-regular if, and only if, $R$ is $*$-abelian, $J(R)$ is nil and $R/J(R)$ is strongly $\pi$-$*$-regular.
\end{proposition}

\begin{proof} The ``only if" part follows like this: Since $R$ is certainly strongly $\pi$-regular, the ideal $J(R)$ is known to be nil (see, e.g., \cite{L}). That the ring $R$ is $*$-abelian follows automatically from Theorem~\ref{5}(3). The final conclusion that $R/J(R)$ retains the same property as that of $R$ follows immediately from Lemma~\ref{7} and Theorem~\ref{5}(2).

Now, we show the validity of the ``if" part. To that aim, since the quotient-ring $R/J(R)$ is strongly $\pi$-$*$-regular, all idempotents of $R/J(R)$ are projections. So, for any, $e\in Id(R)$ it must be that $e+J(R)$ is a projection of $R/J(R)$. Now, since $J(R)$ is nil, with Lemma~\ref{7} at hand there is a projection $p \in R$ such that $e - p \in J(R)$. Therefore, $e(1 - p) \in J(R)$ and $p(1 - e) \in J(R)$ as $R$ is $*$-abelian. It follows at once that $e = ep = pe = p$. Thus, every idempotent of $R$ is a projection. But since $R/J(R)$ is strongly $\pi$-regular, it follows from \cite[Theorem 3]{B} that $R$ is itself $\pi$-regular. Consequently, in view of Theorem~\ref{5}, the ring $R$ is strongly $\pi$-$*$-regular, as claimed.
\end{proof}

The next comments shed some more light on the established above assertions.

\begin{remark}\label{9} 
In regard to the introductory section, it is of some interest to extract some relationships between the different classes of $*$-like rings. Specifically, with the aid of Theorem~\ref{5}, Proposition~\ref{8} and the corresponding results from \cite{LZ}, one infers with no difficulty that strongly $\pi$-$*$-regular rings are themselves strongly $*$-clean which is, in fact, the ``star analogy" with the classical assertion that strongly $\pi$-regular rings are always strongly clean (see, for example, \cite[Theorem 1]{N}). Even something more, it follows from \cite[Corollary 3.8]{CW} that a $*$-ring is strongly $\pi$-$*$-regular if, and only if, it is strongly $*$-clean and $\pi$-regular as well as this is somewhat happened in the case of arbitrary rings (i.e., a ring is strongly clean and $\pi$-regular if, and only if, it is strongly $\pi$-regular).

Moreover, in conjunction with \cite[Theorem 4.9]{CD}, it follows at once that $*$-periodic rings are necessarily strongly $\pi$-$*$-regular which is analogously to the claim that periodic rings are strongly $\pi$-regular. That is why, it could be happen that, for any $*$-ring, the property of being simultaneously $*$-clean and $\pi$-$*$-regular will, eventually, lead to the conclusion that the $*$-ring is strongly $\pi$-$*$-regular.
\end{remark}

We end our work with the following two still unanswered questions.

It is a well-known principal fact that if $R$ is a commutative $\pi$-regular ring, then the $n\times n$ full matrix ring $\mathbb{M}_n(R)$ is strongly $\pi$-regular. So, we come to the following.

\begin{problem}\label{10} Suppose $n\in \mathbb{N}$ and $R$ is a $*$-ring. Determine when $\mathbb{M}_n(R)$ is strongly $\pi$-$*$-regular.
\end{problem}

\begin{problem}\label{11} Characterize \emph{$*$-unit-regular} rings as being those $*$-rings $R$ for which, for each $a\in R$, $a=p+u$ and $aR\cap pR=\{0\}$, where $p^2=p=p^*$ and $u\in U(R)$.
\end{problem}

\section*{Acknowledgments}
The authors are deeply grateful to the anonymous referee for his/her numerous insightful suggestions leading to an improvement of the shape of the article. 

The work of the
first named author J.~Cui is supported by Anhui Provincial Natural
Science Foundation (Grant No. 2008085MA06) and the Key project of
Anhui Education Committee (Grant No. gxyqZD2019009)

 The work of the second named author P.V.~Danchev is partially supported by the Bulgarian National Science Fund under Grant KP-06 No 32/1 of December 07, 2019.

\EditInfo{November 28, 2020}{February 22, 2021}{Pasha Zusmanovich}

\end{paper}